\documentclass[11pt,twoside,letter]{article}
\usepackage[dvips]{graphicx}
\usepackage{amsmath}
\usepackage{amsthm}
\usepackage{amssymb,amscd}
\usepackage{fancyhdr}

\theoremstyle{plain}
\newtheorem{theorem}{Theorem}
\newtheorem{lemma}{Lemma}
\newtheorem{proposition}{Proposition}

\newtheorem{example}{Example}
\newcommand{\R}{{\ensuremath{\mathbb{R}}}}
\newcommand{\Z}{{\ensuremath{\mathbb{Z}}}}
\newcommand{\N}{{\ensuremath{\mathbb{N}}}}

\theoremstyle{definition}
\newtheorem{definition}{Definition}
\theoremstyle{observation}

\newtheorem*{notation}{Notation}

\newtheorem{remark}{Remark}

\begin{document}

\title{\bf\vspace{-39pt} Functions in Sampling Spaces}


\author{Cristina Blanco \\ \small Depto. de Matem\'atica, Univ. de Buenos Aires\\ \small Cdad.
Univ., Pab. I\\ \small  1428 Capital Federal, Argentina \\ \small
cblanco@dm.uba.ar \\
\\
Carlos Cabrelli \\ \small Depto. de Matem\'atica, Univ. de Buenos
Aires and CONICET\\ \small Cdad. Univ., Pab. I \\\small  1428 Capital
Federal,
Argentina \\ \small cabrelli@dm.uba.ar \\
\\
Sigrid Heineken \\ \small Depto. de Matem\'atica, Univ. de Buenos
Aires and CONICET \\ \small  Cdad. Univ., Pab. I \\\small 1428 Capital
Federal, Argentina
\\ \small sheinek@dm.uba.ar }


\date{\today}
\maketitle

\thispagestyle{fancy}

\markboth{\footnotesize \rm \hfill C. BLANCO, C.CABRELLI AND S. HEINEKEN
\hfill} {\footnotesize \rm \hfill FUNCTIONS IN SAMPLING
SPACES \hfill}

\begin{abstract}

Sampling theory in spaces other than the space of band-limited
functions has recently received considerable attention. This is in
part because the  band-limitedness assumption is not very
realistic in many applications. In addition, band-limited
functions can have very slow decay which translates in poor
reconstruction. In this article we study the sampling problem in
general shift invariant spaces. We characterize the functions in
these spaces and provide necessary and sufficient conditions for a
function in $L^2(\R)$ to belong to a sampling space. Furthermore
we obtain decompositions of a sampling space in sampling
subspaces. These decompositions are related with determining sets.
Some examples are provided.

\vspace{5mm}

\noindent{\it Key words and phrases}: Shift invariant spaces, sampling spaces,
frames, determining sets

\vspace{3mm}

\noindent{\it 2000 AMS Mathematics Subject Classification} --- 39A10, 42C40,
41A15
\end{abstract}

\section{Introduction}

Sampling theory originally has been developed in the space
$P_{\sigma}$ of $\sigma$ band-limited functions, i.e. the space of
functions $f\in L^2(\R)$ such that $supp\,\hat{f} \subset
[-\sigma,\sigma],$ where $\sigma > 0$ and
$$\hat{f}(\omega)=\int_{-\infty}^{+\infty}f(x)e^{-2\pi
ix\omega}\,dx.$$
The Classical Sampling
 theorem states that every function $f$ in this class
can be recovered from its samples
${f(\frac{k}{2\sigma})}_{k\in\Z}$. Furthermore the following
reconstruction formula holds:
$$f(x)=\sum_{k\in\Z}f(\frac{k}{2\sigma})
\frac{2\sigma\sin(\pi(2\sigma x-k))}{\pi(2\sigma x-k)},$$ where
the series on the right hand side
 converges uniformly and in $L^2(\R)$.
This theorem has been generalized in many ways. For example to
answer the question whether a similar result can be obtained for
the case of irregular sampling (i.e. when the samples values are
not on a lattice), or when we only have average sampling values,
(i.e. the function has to be recovered from the values
$<f,\phi_k>,$ where each $\phi_k$ is a $L^2$-function that is well
localized around $k$). There is a profuse literature on the
subject. We recommend the books \cite{BF01},\cite{BZ04} and the
survey \cite{AG01} and the references there in.

The space $P_{1/2}$ where the Classical Sampling Theorem holds, is
in particular a Shift Invariant Space (SIS) with generating
function $\psi(x)=\frac{\sin\pi x}{\pi x},$
 moreover $\{\psi(\cdot-k)\}_{k\in\Z}$
is an orthonormal basis of $P_{1/2},$
in particular a frame for $P_{1/2}.$\\
Let us recall that a SIS in $L^2(\R)$ is a subspace  of $L^2(\R)$,
such that it is invariant under integer translations.

Given functions $f_1,f_2,\dots ,f_n \in L^2(\R)$, we will denote by
$S(f_1,\dots,f_n)$, the SIS generated
by the integer translates of these functions, i.e. the $L^2(\R)$-
closure of the span of the set $\{f_i(\cdot-k): i=1,...,n, k \in \Z\}$.

When the  set $\{f_i(\cdot-k): i=1,...,n, k \in \Z\}$,\ forms a
frame of $S(f_1,\dots,f_n)$, we will write sometimes
$V(f_1,\dots,f_n)$ instead of $S(f_1,\dots,f_n)$, to stress this
fact. It is known that every space $S(f_1,\dots,f_n)$ contains
functions $g_1,\dots,g_l, $
 with $l \leq n,$ such that
 $S(f_1,\dots,f_n) = V(g_1,\dots,g_l)$.

There is a natural interest in trying to obtain sampling theorems
in SIS's other than spaces of band-limited functions,
\cite{Wal92},\cite{SZ99}, since the band-limitedness assumption is
not realistic in  numerous applications in digital signal and
image processing, \cite{AG01}.

Our main interest in this article is to study the structure
and find characterizations of functions of $L^2(\R)$
that belong to a general sampling space $V(\varphi)$.
A sampling space is a space $V(\varphi)$ where the generator
has special properties.

Throughout this paper we will consider  the following
class of sampling spaces:
\begin{definition}

A closed subspace $V(\varphi)$ of $L^2(\R)$ is called a {\sl sampling
space} if there exists a function $s$ such that:

\begin{enumerate}
\item The translates $\{s(\cdot-k)\}_{k\in\Z}$ are a frame for the
space $V(\varphi).$

\item For every sequence $\{c_k\}_{k\in \Z}\in \ell^2(\Z)$ the
series $\sum_{k\in \Z} c_k s(\cdot-k)$ converges pointwise to a
continuous function.

\item \label{continuity} For every $f\in V(\varphi),$
\begin{equation}\label{desf}
f(x)=\sum_{k\in \Z} f(k) s(x-k),
\end{equation}
where the convergence is in $L^2(\R)$ and uniform in $\R.$
\end{enumerate}
The function $s$ is called the {\sl sampling function} of
$V(\varphi),$ and is the only function in $V(\varphi)$ that
satisfies equation (\ref{desf}). (See Lemma \ref{L3}).
\end{definition}

We found that if a function $f$ belongs to a sampling space, then
$S(f)$ is a sampling space itself. We characterize the sampling
function of $S(f)$ as the orthogonal projection, onto $S(f)$, of
the sampling function of the original space. From here, using a
characterization of  sampling spaces given by Sun and Zhou
\cite{SZ99}, we obtain necessary and sufficient conditions for a
function $f$ to belong to a sampling space.

These results can be connected with the problem of {\sl
Determining Sets} for sampling spaces, \cite{ACHMR03}. Basically
a set of functions of a finitely generated SIS is a determining
set if and only if it is a set of generators (in the sense
mentioned above) of the SIS. Necessary and sufficient conditions
for a finite subset of a general finitely generated SIS to be a
determining set have been obtained in \cite{ACHMR03}. We obtain
the equivalent conditions in the case of sampling spaces and we
show that given a determining set of a sampling space, then the
space can be decomposed as a sum of smaller sampling spaces, each
of them generated by one of the functions in the determining set.
We characterize all possible (countable) decompositions of a
sampling space $V(\varphi)$ as a direct sum of sampling spaces.

The grammian of a function $\varphi \in L^2(\R)$ is the function
$G_{\varphi}(\omega) = \sum_{k \in \Z}|\hat{\varphi}(\omega +
k)|^2.$ Denote by $E_\varphi$ the set $E_\varphi= \{ w \in
\R:G_\varphi(\omega)
>0 \}.$ The set $E_\varphi$ is periodic i.e. $E_\varphi=E_\varphi+k$ for
every integer $k.$ If $\psi$ is any generator of $S(\varphi),$
then $E_\psi=E_\varphi.$ If $V(\varphi)$ is a sampling space,
every partition of $E_\varphi$ into countable measurable periodic
sets defines a decomposition of $V(\varphi)$ as a direct sum of
sampling spaces. The sampling function of each of these spaces
turns out to be the orthogonal projection of the sampling function
over each of the components. This type of decompositions has been
obtained in \cite{BDR94},\cite{Bow00}, for SIS of $L^2(\R)$. We
show here that it is still valid for sampling spaces, that is, the
subspaces in the decomposition are also sampling spaces.

In \cite{SZ04}, Sun and Zhou obtained sufficient conditions for  a function
to belong to a sampling space. The conditions are not necessary as
we show in Example~\ref{Ej2}. We found necessary and sufficient
conditions for a function to belong to a sampling space, for the
case of sampling spaces with generators having Fourier transforms
in $L^1(\R)$. As expected, the $\sigma$ band-limited functions
satisfy these conditions (see Example~\ref{Ej1}).

The paper is organized as follows: In section 2 we state the
 main results, in section 3 we present the proofs and section 4 is
devoted to some examples.

\section{Statement of Results}

\begin{definition}
A sequence $\{\psi_k\}_{k\in\Z}$ is a {\sl frame} for a separable
Hilbert space $\mathcal{H}$ if there exist positive constants $A$
and $B$ that satisfy
$$A\|f\|^2\leq\sum_{k\in\Z}|\langle f,\psi_k\rangle|^2 \leq
B\|f\|^2\,\,\,\,\,\forall f \in \mathcal{H}.$$

If $\{\psi_k\}_{k\in\Z}$ satisfies the right inequality in the
above formula, it is called a {\sl Bessel sequence}.
\end{definition}

For a shift invariant space $S(\psi)$ it is known (\cite{BDR94})
that the integer translates of the function $\varphi$ defined by

\begin{equation}\label{g}
\hat{\varphi}(\omega)=
\begin{cases}
\frac{\hat{\psi}(\omega )}{G_{\psi}(\omega)^{\frac{1}{2}}}&
\text{for } \omega \in E_{\psi} \\
 0 & \text{otherwise},
\end{cases}
\end{equation}
form a tight frame of $S(\psi),$ in particular
$S(\psi)=V(\varphi).$
Now we state our first theorem:
\begin{theorem}\label{teorsf}
Let $V(\varphi)$ be a sampling space with sampling function $s$
and $f\in L^2(\R).$ If $f\in V(\varphi)$ then $S(f)$ is a sampling
space with  sampling function $s_f=P_{S(f)}(s),$ where
$P_{S(f)}(s)$ is the orthogonal projection of $s$ onto $S(f).$ In
this case $\hat{s_f}=\hat{s} \chi_{E_f}.$
\end{theorem}
We will use this result, together with the characterization of a
sampling space obtained by Sun and Zhou in \cite{SZ99}, to obtain
necessary and sufficient conditions on a function $f$ to belong to
a sampling space.

We first need the following definition:
\begin{definition}
For $f\in L^2(\R)$ the Zak transform of $f$ is the function on $\R^2$:
$$Z_f(x, \omega)=\sum_{k\in \Z}f(x+k) e^{-2\pi i k \omega}.$$
\end{definition}
For properties of the Zak transform see \cite{Jan88}.

\begin{theorem}\label{T2}
Assume $f \in L^2(\R).$ Then $f$ belongs to a sampling space if
and only if the function $h$ defined by $\hat{h}=\hat{f}/ G_f$ in
$E_f$ and zero otherwise, satisfies:
\begin{enumerate}

\item $h$ is continuous.

\item The function$\sum_{k\in\Z}|h(x-k)|^2$ is bounded on $\R.$
\item There exist constants $A, B>0$ such that
$$A\chi_{E_{f}}(\omega) \leq |Z_{h}(0, \omega)| \leq
B\chi_{E_{f}}(\omega)\,\,\,\,\,\text{a.e. } \omega.
$$

\end{enumerate}
\end{theorem}

This results can now be applied to the problem of {\sl determining
sets}.

Given a set $\mathcal{F}=\{f_1,\dots,f_m\}\subset L^2(\R)$ where
$f_i$ are functions that belong to an unknown shift invariant
space $V(\varphi),$ it is an important matter to be able to decide
whether this set $\mathcal{F}$ is sufficient to determine
$V(\varphi).$ This leads to the concept of {\sl determining sets}
of shift invariant spaces. In \cite{ACHMR03}, the problem is
solved if $\{\varphi(\cdot-k)\}_{k\in\Z}$ is a Riesz basis of
$V(\varphi).$

 We study the problem in the case that $V(\varphi)$ is a sampling
  space. We give necessary and sufficient conditions on $\mathcal{F},$
needed for determining the unknown sampling space $V(\varphi).$
Moreover we decompose the sampling space $V(\varphi)$ as the sum
of the sampling spaces $S(f_i).$ In particular the sampling
function $s$ of $V(\varphi)$ can be recovered from the functions
$f_i.$

\begin{definition}
Let $V(\varphi)$ be  a shift invariant space. The set
$\mathcal{F}=\{f_1,\dots,f_m\} \subset V(\varphi)$ is a {\sl
determining set} for $V(\varphi)$ if for any $g\in V(\varphi)$
there exist $\hat{\alpha_1}, \dots, \hat{\alpha_m}$ 1-periodic
measurable functions such that:
$$\hat{g}=\hat{\alpha_1}\hat{f_1}+\cdots +\hat{\alpha_m}\hat{f_m}.$$
\end{definition}

\begin{theorem}\label{pro}
Assume $V(\varphi)$ is a sampling space. Then
$\mathcal{F}=\{f_1,\dots,f_m\}\subset V(\varphi)$ is a determining
set for $V(\varphi)$ if and only if the set
$$Z=\left(\bigcup_{i=1}^m E_{f_i}\right) \triangle E_{\varphi}$$ has Lebesgue
measure zero (where $\triangle$ denotes the symmetric difference
of sets).

Moreover, if $\mathcal{F}$ is a determining set for $V(\varphi)$,
then
$$V(\varphi)=S(f_1)+ \dots + S(f_m).$$
\end{theorem}

From Theorem \ref{teorsf} we can also obtain the following
decomposition of a sampling space.
\begin{proposition}\label{prop1}
Let  $V(\varphi)$ be a sampling space
 and let $\{E_j\}_{j \in N}$ be a partition of $E_{\varphi}$ in
 periodic measurable sets (i.e. $E_{\varphi} = \bigcup_j E_j$,
 $ |E_j \cap E_l| =0 $ for all $j \neq l,$ and $E_j+k=E_j$ for
 every integer $k).$

Define for $j \in \N,$ $\varphi_j$ by $\hat\varphi_j = \hat
\varphi \chi_{E_j}$. Then we have:

$$V(\varphi) = \bigoplus_j V(\varphi_j),$$  where $V(\varphi_j)$ is a
sampling space for each $j$. Furthermore, if $s$ is the sampling
function of $V(\varphi)$ then the Fourier transform of the
sampling function $s_j$ of $V(\varphi_j)$ is $\hat s_j=\hat
s\chi_{E_j}.$
\end{proposition}

We now return to the study of the characterization of functions
that belong to sampling spaces. Our results in this part are
related to a recent result of Sun and Zhou \cite{SZ04}. Let
$$ \mathcal{U}=\{f\in L^2(\R): \hat{f}\in L^1(\R)\}. $$
In \cite{SZ04}, sufficient conditions for a function $f$ in
$\mathcal{U}$ to belong to a sampling space are given. We show in
Example~\ref{Ej2} that these conditions are not necessary.

We found necessary and sufficient conditions for a function in
$\mathcal{U}$ to belong to a sampling space for the case that the
generator of the sampling space is in $\mathcal{U}.$ Note that if
the generator $\varphi$  of the sampling space $V(\varphi)$ is in
$\mathcal{U},$ then any other function $\psi$ such that
$V(\psi)=V(\varphi)$ is also in $\mathcal{U}.$

The conditions are stated in the following theorem:

\begin{theorem}\label{t1}

 Let $f\in \mathcal {U}.$ Then  $f$ belongs to a sampling space
 $V(\varphi)$
 with $\varphi\in\mathcal{U}$
 if and only if $f$ satisfies:

\begin{enumerate}
\item[a)]The sequence $\{f(k)\}_{k\in\Z} \in \ell^2(\Z)$,

\item[b)]there exist $A, B >0$ such that
$$ A |Z_f(0,\omega)|^2 \leq G_f(\omega)\leq B |Z_f(0,\omega)|^2
 \;\; \text{ a.e.}
\; \omega \in \R,$$

\item[c)]The integral $$\int_{[0,1]\cap
E_f}\frac{\sum_{k\in\Z}|\hat{f}(\omega+k)|} {|Z_f(0,\omega)|}
d\omega <+\infty,$$

\item[d)]There exists $L > 0$ such that
$$\int_{[0,1]\cap E_f}\left|\frac{Z_{\hat{f}}(\omega, -x)}{Z_f(0,\omega)}\right|^2 d\omega<L
< + \infty\,\,\,\forall x \in \R.$$

\end{enumerate}
\end{theorem}

\begin{remark}
We have, in fact, the following stronger version of the necessary
condition in Theorem~\ref{t1}:
\end{remark}

\begin{proposition}\label{nece}
If $f\in L^2(\R)$ belongs to a sampling space $V(\varphi)$ with
$\varphi\in\mathcal{U},$ then $ f$ satisfies $a),\, b),\, c),\,
d).$
\end{proposition}

Finally we observe that it is natural to consider the sampling
problem on a general lattice $a\Z+b$, with $a\in\R_{>0}, b \in
\R.$ For this, let $T:L^2(\R) \rightarrow L^2(\R)$ be an unitary
operator and $V=V(\varphi)$ a sampling space. Set $W=T(V)$.

Consider $g \in W$ and call $f = T^{-1}g$. Then $f=\sum_k f(k)
t_k\varphi, $ where $t_k$ is the translation operator
$t_kh(x)=h(x-k)$. So, $g=Tf=\sum_kf(k)(T\circ t_k)\varphi=\sum_k
T^{-1}(Tf)(k) (T\circ t_k)\varphi.$ That is
$$ g= \sum_k (T^{-1}g)(k) (T \circ t_k) \varphi, \;\;\; \forall g \in W. $$

Let us now define the unitary dilation operator $D_a$ by
$D_af(x)=\sqrt{a}f(ax)$ and $T$ by $T=D_a \circ t_b. $ Denote
$\psi = T\varphi. $ Then, due to the commutation relation $ D_a
t_k = t_{\frac{k}{a}} D_a, $ the following sampling formula holds:

$$ g= \sum_k g(\frac{k+b}{a}) \psi(x- \frac {k} {a}), \;\;\;\;\;\;
\forall g \in V_{a,b} = \overline{span}(\{\psi(\cdot-\frac{k}{a}), k\in
\Z\})). $$

\section{Proofs}
We will first state some known results that we will use in the
proofs.

In \cite{BL22} Benedetto and Li proved the following result:
\begin{proposition}
The sequence $\{\psi(\cdot-k)\}_{k\in\Z}$ is a frame for the
closure in $L^2(\R)$ of the space it spans if and only if there
exist positive constants $A$ and $B$ that satisfy
$$A \leq G_{\psi}(\omega)\leq B \,\,\,\,\,\,\text{ a.e. } \omega \in E_{\psi}.$$
\end{proposition}

\begin{proposition}\label{Ron}(Theorem 2.14 \cite{BDR94})
Let $S(\psi)$ be a shift invariant space. A function $f$ is in
$S(\psi)$ if and only if $\hat{f}=r \hat{\psi}$ for some periodic
function $r$ of period one, with $r \hat{\psi} \in L^2(\R).$
\end{proposition}

We will also need this known result: if the  sequence
$\{\psi_k\}_{k\in\Z}$ is a frame for the Hilbert space
$\mathcal{H},$ then $f\in\mathcal{H} \text{ if and only if }
f=\sum_{k\in \Z}c_k\psi_k, \text{ for some } c_k\in \ell^2(\Z)$
with convergence in $L^2(\R).$

Sun and Zhou gave the following characterization, which we will
use later to determine if a function generates a sampling space:
\begin{proposition}(\cite{SZ99}, Theorem 1)\label{SZ99}
Let $V(\varphi)$ be a shift invariant space. Then the following
two assertions are equivalent:
\begin{enumerate}
\item[(i)]The space $V(\varphi) $ is a sampling space
\bigskip
\item[(ii)] The function $\varphi $ is continuous, $\sum_{k\in \Z}|\varphi(x-k)|^2$ is
bounded on $\R$ and
$$A\chi_{E_{\varphi}}(\omega) \leq |Z_{\varphi}(0, \omega)|
\leq B\chi_{E_{\varphi}}(\omega)\,\,\,\,\,\text{a.e. } \omega$$ for
some constants $A, B>0.$
\end{enumerate}
\end{proposition}

\begin{proposition}\label{Riesz}(\cite{SZ04}, Theorem 1.2)
Let $V(\varphi)$ be a sampling space with sampling function $s.$
Then there exists a sampling space $V(\tilde{\varphi})$ with
sampling function $\tilde{s}$ such that:

\begin{enumerate}
\item The space  $V(\varphi)$ is a subspace of
$V(\tilde{\varphi})$

\item The sequence$\{\tilde{s}(\cdot-k)\}_{k\in\Z}$ is a Riesz
basis for $V(\tilde{\varphi}).$
 \end{enumerate}

 \end{proposition}

\subsection{Proofs of Theorems 1 and 2}

First we need the following lemmas:

\begin{lemma}\label{L1}
Assume $V(\varphi)\subset L^2(\R)$ is a shift invariant space and
$\varphi$ is a continuous function. Let $\psi\in V(\varphi)$ such
that $\{\psi(\cdot-k)\}_{k\in \Z}$ is a Bessel sequence. If
$\sum_{k\in\Z}|\varphi(x+k)|^2<L<+\infty \,\,\,\forall x\in \R,$
then $\sum_{k\in\Z}|\psi(x+k)|^2<L'<+\infty \,\,\,\forall x\in
\R.$
\end{lemma}
The proof of this lemma is in \cite{SZ99} for the case that
$\{\psi(\cdot-k)\}_{k\in \Z}$ is a frame of $V(\varphi),$ but the
same proof works if it is only a Bessel sequence.

\begin{lemma}\label{L2}(\cite{Sun05}, Lemma 2.5)
If $\varphi\in L^2(\R)$ is continuous and
$\sum_{k\in\Z}|\varphi(x+k)|^2<L<+\infty,$ then
$Z_{\varphi}(x,\omega)=0\,\,\,\forall x\in \R ,\,\, \text{ a.e.}
\,\,\omega \in \R\setminus E_{\varphi}.$
\end{lemma}

\begin{remark}\label{cuenta}
Let $V(\varphi)$ be a sampling space and $s$ its sampling
function. For $f\in V(\varphi)$,
$$\hat{f}(\omega)=Z_f(0,\omega)\hat{s}(\omega)\,\,\,\text{
a.e.} \, \omega \in \R.$$
Therefore we obtain
$$G_f(\omega)=|Z_f(0,\omega)|^2G_s(\omega)\,\,\,\text{
a.e.} \, \omega \in \R.$$
\end{remark}

To see this, observe that
$$ f(x) = \sum_{k \in \Z} f(k)s(x-k),$$ with uniform convergence and in
$L^2(\R),$ so

\begin{equation*}
\hat {f} (\omega) = (\sum_{k \in \Z} f(k) e^{-2 \pi i k \omega})
\hat{s} (\omega)=Z_f(0,\omega) \hat{s} (\omega)  \,\,\,\, \text{
a.e.} \, \omega \in \R,
\end{equation*}

and then

$$\sum_k |\hat {f} (\omega+k)|^2=
|Z_f(0,\omega)|^2 \sum_k| \hat{s} (\omega+k)|^2\,\,\,\,\text{
a.e.} \, \omega \in \R,$$ hence
$$G_f(\omega)= |Z_f(0,\omega)|^2  G_{s}(\omega) \,\,\,\,\text{
a.e.} \, \omega \in \R.$$

\begin{lemma}\label{L3}
Let $V(\varphi)$ be a sampling space and $s$ its sampling
function. Then we have:
\begin{enumerate}

\item [i)] If $\psi_1,\psi_2 \in L^2(\R)$ and
$S(\psi_1)=S(\psi_2),$ then $E_{\psi_1}=E_{\psi_2}$ (up to a set
of measure zero). In particular if $\{\psi(\cdot-k)\}_{k \in \Z}$
is a frame for $ V(\varphi)$ then $E_{\varphi}=E_{\psi}.$

\item[ii)] The sampling function satisfies
$Z_s(0,\omega)=\chi_{E_s}(\omega)\,\,\,\, \text{ a.e.} \, \omega
\in \R.$

\item [iii)] For $f\in V(\varphi),\, E_f=\{\omega \in \R:
Z_f(0,\omega)\neq 0\}$ (up to a set of measure zero).

\item [iv)]  The sampling function $s$ is unique, up to a set of
measure zero, and satisfies $\hat s(\omega) =
\frac{\psi(\omega)}{Z_{\psi}(0,\omega)}
\chi_{E_{\varphi}}(\omega)$ for each generator $\psi$ whose
translates form a frame of $V(\varphi)$.

\end{enumerate}
\end{lemma}

\begin{proof}
$i)$ If $\{\psi(\cdot-k)\}_{k \in \Z}$ is a frame for $
V(\varphi)$ then there exists $\{c_k\}_{k\in \Z} \in \ell^2(\Z)$
such that
$$\varphi(x)=\sum_{k\in \Z}c_k\psi(x-k).$$
So $G_{\varphi}(\omega)=|(\sum_{k\in \Z} c_k e^{-2 \pi i  k
\omega})|^2 G_{\psi}(\omega)\,\,\text{ a.e. } \omega \in \R,$
which yields $E_\varphi\subseteq E_\psi$ (up to a set of measure
zero). Similarly we obtain the other inclusion. Now, if $\psi_1,
\psi_2$ generate the same SIS, we can modify these generators as
in (\ref{g}) to obtain tight frames, and the result follows.

$ii)$ Since we have
$\hat{s}(\omega+k)=Z_s(0,\omega)\hat{s}(\omega+k)
 \,\,\text{ a.e. } \omega \in \R, $ it follows that
$Z_s(0,\omega)=1$ for almost every $\omega \in E_s.$ On the other
side, by Proposition~\ref{SZ99}, $Z_s(0,\omega)=0$ for almost
every $\omega \notin E_s.$

$iii)$ We have
$$G_f(\omega)= |Z_f(0,\omega)|^2  G_{s}(\omega) \,\,\,\,\text{
a.e.} \, \omega \in \R.$$ Using Proposition~\ref{Riesz} we can
assume that $\{s(\cdot-k)\}_{k\in\Z}$ is a Riesz basis of the
sampling space, so $G_s(\omega) \neq 0 \text{ a.e. }\,\omega \in
\R,$ which implies $iii).$

$iv)$ Let $\psi$ be a generator whose translates form a frame of
$V(\varphi).$ We have $\hat \psi (\omega) = Z_{\psi}(0,\omega)
\hat s(\omega) \,\,\text{ a.e. } \omega \in \R.$ By $i),\,
E_{\varphi}=E_{\psi},$ so using $iii),$

$$\hat s (\omega) = \frac{\hat
\psi (\omega)}{ Z_{\psi}(0,\omega)}\,\,\text{ a.e. } \omega \in
E_{\varphi}$$
 and the result follows.

\end{proof}

 \begin{notation}
 Let $f,g\in L^2(\R).$ We denote
 $$[f,g](x)=\sum_{k\in\Z}f(x+k)\overline{g(x+k)}.$$
 \end{notation}
 Observe that $G_f(\omega)=[\hat{f},\hat{f}](\omega).$

Since the Fourier transform preserves the scalar product, if $V$
is a closed subspace of $L^2(\R)$ and $P_V$ is its orthogonal
projection, we have that $\widehat{P_V(f)}=P_{\hat{V}}(\hat{f}),$
where $\hat{V}=\{\hat{f}:f\in V\}.$

In \cite{BDR94} the following formula for the orthogonal
projection was obtained:

$$P_{\widehat{S(\varphi)}}(\hat{f})(\omega)=r(\omega)\hat{\varphi}(\omega),$$
where
\begin{equation}
r(\omega)=
\begin{cases}
\frac{[\hat{f},\hat{\varphi}](\omega)}
{[\hat{\varphi},\hat{\varphi}](\omega)}&
\text{for } \omega \in E_\varphi \\
 0 & \text{otherwise. }
\end{cases}
\end{equation}

We are finally ready to prove the theorems.
\begin{proof}[Proof of Theorem 1]
 For $f\in V(\varphi),$ we define
\begin{equation}\label{h}
\hat{h}(\omega)=
\begin{cases}
\frac{\hat{f}(\omega )}{Z_f(0,\omega)}&
\text{for } \omega \in E_f \\
 0 & \text{otherwise. }
\end{cases}
\end{equation}
Observe that by Lemma~\ref{L3} $iv)$ the function $h$ is well
defined. Using Proposition~\ref{Ron}, it is easy to see  that
$S(f)=S(h).$ To show that $\{h(\cdot - k )\}_{k \in \Z}$ is a
frame sequence, first observe that $E_h=E_f,$ so for almost every
$\omega \in E_h$,
$$ G_h(\omega)=\frac{G_f(\omega)}{|Z_f(0,\omega)|^2}=
G_s(\omega),$$ which is uniformly bounded below and above.

Now we will see that $S(f)=V(h)$ is a sampling space:

Since $h \in S(f) \subseteq V(\varphi),$ $h$ is continuous. By
Proposition~\ref{SZ99} and Lemma~\ref{L1}, $\sum_{k \in \Z} |h(x
+k)|^2$ is uniformly bounded in $\R.$ So by Lemma~\ref{L2},
$Z_h(0,\omega)=0 \,\,\, \text{ a.e.}\, \omega \in \R\setminus E_h
.$ Furthermore, for almost every $\omega \in E_h$,
$$G_h(\omega)=|Z_h(0,\omega)|^2 G_s(\omega).$$
Since $E_h \subseteq E_s$ (up to a set of measure zero), we can
write for almost all $\omega \in E_h$
$$\frac{G_h(\omega)}{G_s(\omega)}= |Z_h(0,\omega)|^2,$$
and using that $G_h(\omega)$ and $G_s(\omega)$ are both bounded
above and below in $E_h,$ we have that $|Z_h(0,\omega)|$ is also
bounded above and below. Hence, using Proposition~\ref{SZ99}, we
can conclude that $S(f)$ is a sampling space. Moreover, $h$ is its
sampling function. To see this,
it suffices to prove that for every $g\in
S(f),\, \hat{g}(\omega)=Z_g(0,\omega) \hat{h}(\omega).$ But for $g
\in S(f)\subseteq V(\varphi),$
$$\hat{g}(\omega)=Z_g(0,\omega)\hat{s}(\omega),$$
and since $\hat{h}(\omega)=\hat{s}(\omega)$ for $\omega\in
E_{\varphi},$ we have that $\hat{g}(\omega)=Z_g(0,\omega)
\hat{h}(\omega)$ for almost every $\omega \in E_f.$ On the other
hand, since $g\in S(f),$ there exists a 1-periodic function $r,$
such that $\hat{g}(\omega)=r(\omega) \hat{f}(\omega),$ and
therefore the equality also holds for almost every $\omega \in \R
\setminus E_f.$

Hence
$$\hat{s_f}(\omega)=\hat{h}(\omega)=
\hat{s}(\omega)\chi_{E_f}(\omega) \,\,\text{ a.e. }\,\omega \in
\R.$$

Finally we prove that $\hat{h}$ is actually the projection of $s$
onto $S(f),$ i.e. $\hat{h}=\widehat{P_{S(f)}(s)}.$

For almost every $\omega \in E_f$,
$$\widehat{P_{S(f)}(s)}(\omega)=\frac{[\hat{s},\hat{f}](\omega)}{[\hat{f},\hat{f}](\omega)}
\hat{f}(\omega)=$$\\
$$\frac{[\hat{s},Z_f(0,\cdot)\hat{s}](\omega)}{[\hat{f},\hat{f}](\omega)}
Z_f(0,\omega)\hat{h}(\omega)=\\
|Z_f(0,\omega)|^2
\frac{[\hat{s},\hat{s}](\omega)}{[\hat{f},\hat{f}](\omega)}
\hat{h}(\omega)=\hat{h}(\omega),$$ and for almost every $
\omega\in \R\setminus E_f,$ we have
$$\widehat{P_{S(f)}(s)}(\omega)=\hat{h}(\omega)=0.$$
This completes the proof of Theorem~\ref{teorsf}.
\end{proof}
Theorem~\ref{T2} is an immediate consequence of
Proposition~\ref{SZ99}.

\subsection{Proofs of Theorem \ref{pro} and Proposition \ref{prop1}}

\bigskip

\begin{proof}[Proof of Theorem ~\ref{pro}]Let $s$ be the sampling function of $V(\varphi).$\\
Assume that $\mathcal{F}=\{f_1,\dots,f_m\}\subset V(\varphi)$ is a
determining set for $V(\varphi).$ Recall that by Lemma~\ref{L3},
$E_{\varphi}=E_s$ (up to a set of measure zero). Since
$$G_{f_i}(\omega)=|Z_{f_i}(0,\omega)|^2G_s(\omega) \,\,\, \text{ a.e.
} \omega \in \R,$$ for almost every $\omega \in E_{f_i},\, \omega$
belongs to $E_{\varphi},$ so the set $\cup_{i=1}^m E_{f_i}
\setminus E_{\varphi}$ has Lebesgue measure zero.

On the other side there exist $\hat{\alpha_1}, \dots,
\hat{\alpha_m}$ 1-periodic measurable functions such that
$$\hat{\varphi}(\omega)=\sum_{i=1}^m \hat{\alpha_i}(\omega)\hat{f_i}(\omega),$$
so
$$G_{\varphi}(\omega)=\sum_{k\in \Z} |\hat{\varphi}(\omega+k)|^2 \leq
\sum_{k\in \Z} \left(\sum_{i=1}^m
|\hat{\alpha_i}(\omega)\hat{f_i}(\omega+k)| \right)^2.$$ Hence the
set  $E_{\varphi} \setminus \cup_{i=1}^m E_{f_i}$ has Lebesgue
measure zero.

To prove the reciprocal it suffices to show that there exist
  $\hat{\alpha_1}, \dots,
\hat{\alpha_m}$ 1-periodic measurable functions such that
\begin{equation}\label{E}
\hat{s}(\omega)=\sum_{i=1}^m
\hat{\alpha_i}(\omega)\hat{f_i}(\omega).
\end{equation}
 Define the sets $B_i$
inductively by $B_1= E_{f_1},$ and for $2\leq i \leq m,\,\,\ B_i=
E_{f_i} \setminus \bigcup_{j=1}^{i-1} B_j.$ For $1\leq i \leq m$
set
\begin{equation}
\hat{\alpha_i}(\omega):=
\begin{cases}
\frac{1}{Z_{f_i}(0,\omega)}&
\text{for } \omega \in B_i \\
 0 & \text{otherwise. }
\end{cases}
\end{equation}

Since $\hat{s}(\omega)=\hat{\alpha_i}(\omega)
\hat{f_i}(\omega)\,\,\,\text{ a.e. } \omega \in B_i,$ equation
(\ref{E}) holds.

Finally we will see that $V(\varphi)=S(f_1)+ \dots + S(f_m).$
Since
$$\int_{\R}|\hat{\alpha_i}(\omega)|^2 |\hat{f_i}(\omega)|^2 d\omega=\int_{B_i}
|\hat{s}(\omega)|^2 d\omega < +\infty,$$

it follows that $\hat{\alpha_i} \hat{f_i} \in S(f_i) $ (see
Proposition~\ref{Ron}), so $V(\varphi) \subseteq S(f_1)+\dots +
S(f_m).$

To show the other inclusion, let $g\in S(f_1)+\dots + S(f_m).$
Then there exist $\hat{\beta_1}, \dots, \hat{\beta_m}$ 1-periodic
measurable functions such that $\hat{\beta_i}\hat{f_i}\in
L^2(\R),\,\,1\leq i \leq m ,$ and
$\hat{g}=\hat{\beta_1}\hat{f_1}+\cdots +\hat{\beta_m}\hat{f_m}.$
So, using that
$\hat{f_i}(\omega)=Z_{f_i}(0,\omega)\hat{s}(\omega),$ we have that
$\hat{g}(\omega)=(\hat{\beta_1}(\omega)Z_{f_1}(0,\omega)+\cdots
+\hat{\beta_m}(\omega)Z_{f_m}(0,\omega))\hat{s}(\omega).$ Since
$\hat{g}\in L^2(\R)$ and $\hat{\beta_1}Z_{f_1}(0,\omega)+\cdots
+\hat{\beta_m}Z_{f_m}(0,\omega)$ is a 1-periodic function,
applying again Proposition~\ref{Ron}, $g$ belongs to $V(\varphi).$
\end{proof}

\begin{proof}[ Proof of Proposition ~\ref{prop1}]

Since $\varphi$ satisfies that there exist $A,B \geq 0$ such that
$$A \leq G_{\varphi} (\omega)\leq B \,\,\, \text{ a.e.
} \omega \in E_\varphi ,$$ then each $\varphi_j$ has the same
property in $E_j$. We conclude that $\{\varphi_j(\cdot  -
k)\}_{k\in \Z}$ is a frame sequence. Furthermore, note that
$\varphi_j\in V(\varphi)$ since $\chi_{E_j}$ is 1-periodic and
$\chi_{E_j}\hat{\varphi}\in L^2(\R).$ Hence, by by
Theorem~\ref{teorsf}, $V(\varphi_j)$ is a sampling space.

Let us show now that $V(\varphi) = \bigoplus_j V(\varphi_j).$
Assume that $f \in V(\varphi).$ Then $\hat f = m_f  \hat\varphi$
with $m_f \in L^2[0,1).$

So, $\hat f = \sum_j (m_f \chi_{E_j}) ( \chi_{E_j}
 \hat \varphi) = \sum_j m_j \hat \varphi_j,$
where $m_j = m_f \chi_{E_j} \in L^2[0,1).$ Then $\hat f_j :=
m_j\hat\varphi_j
 \in V(\varphi_j).$ That is $f = \sum_j f_j, f_j \in V(\varphi_j).$

On the other side, assume that $g_j \in V(\varphi_j)$ and $\sum_j
g_j = 0.$ Write $\hat g_j = \theta_j \hat \varphi_j$ with
$\theta_j \in L^2[0,1).$

Suppose that for some $r \in \N,$ the set $ M := \{\omega: \hat
g_r (\omega) \neq 0 \}$ has positive Lebesgue measure. Since $M
\subset E_r$ we have that for almost all $\omega\in M,$

$$ 0 = \sum_j \hat g_j(\omega) = \hat g_r(\omega) =
\theta_r(\omega) \hat \varphi_r(\omega).$$

The fact that $\omega  \in E_r$ implies that for some integer $k$,
$\hat\varphi(\omega+k) \neq 0.$

Then we can write $0 = \theta_r(\omega+k) \hat
\varphi_r(\omega+k)= \theta_r(\omega) \hat \varphi_r(\omega+k).$

So, $\theta_r(\omega)=0.$ That is $\theta_r\equiv 0 \,\,\ \text{
a.e. }$ in $M,$ which is a contradiction. We conclude that $g_j
\equiv 0$ for all $j.$

Since $\varphi_j\in V(\varphi),$ by Theorem~\ref{teorsf},
$\hat{s_j}=\hat{s}\chi_{E_{\varphi_j}}=\hat{s}\chi_{E_j}$ and this
completes the proof.

\end{proof}

\subsection{Proof of Theorem \ref{t1}} First note that if in
Proposition~\ref{Riesz} we suppose that the sampling function $s$
of $V(\varphi)$ belongs to $\mathcal{U},$  it follows from the
proof that the sampling function $\tilde{s}$  of
$V(\tilde{\varphi})$ also belongs to $\mathcal{U}.$

Next we prove the following lemma:

 \begin{lemma}\label{lemma}
 Let $f$ be a  function in $L^2(\R).$
 If $f$ belongs to a sampling space,
 then the sequence $\{f(k)\}_{k \in \Z}$ belongs to $\ell^2(\Z)$.
\end{lemma}

 \begin{proof}
 If $f$ belongs to a sampling space, by Proposition~\ref{Riesz} we can assume that
  $f$ belongs to a sampling space $V(\varphi)$
where the integer translates of the sampling function $s$ is a
Riesz sequence. That is, there exist constants $A,B > 0$ such that
\begin{equation}\label{frame-condition}
A \leq \sum_k| \hat{s} (\omega+k)|^2 \leq B \,\,\,\text{ a.e.} \;
\omega \in \R.
\end{equation}
On the other side, applying Remark~\ref{cuenta}
$$G_f(\omega)=|Z_f(0,\omega)|^2G_s(\omega).$$
Using \eqref{frame-condition},
 the left hand side of the equation
 $$\frac {\sum_k |\hat {f} (\omega+k)|^2} {\sum_k| \hat{s} (\omega+k)|^2} =
  |Z_f(0,\omega)|^2
 \,\,\,\,\text{ a.e.} \,
\omega \in \R,
$$
 has finite integral in $[0,1]$ and coincides
almost everywhere with the  right hand side, so the
function $Z_f(0,w)$ is in $L^2[0,1]$ and then the sequence
$\{f(k)\}_{k \in \Z}$ belongs to $\ell^2(\Z)$.
 \end{proof}

 The following result is a version of the Poisson Summation
 Formula.
\begin{proposition}
If $f\in \mathcal{U}$ and $\{f(k)\}_{k\in \Z} \in \ell^2(\Z)$,
then:

$$Z_f(0,\omega)=\sum_{k\in \Z}\hat{f}(\omega +k)\,\,\,\text{ a.e. } \omega\in [0,1].$$
\end{proposition}
\begin{proof}
Observe that the function $\sum_{k\in \Z}\hat{f}(\omega +k)\in
L^1[0,1]$ and its Fourier coefficients are $\{f(k)\}_{k\in\Z}$ and
so this function is in $L^2[0,1].$
\end{proof}

\begin{remark}
If $V(\varphi)$ is a sampling space with sampling function $s,$
then it follows from  Lemma~\ref{lemma} that
$$V(\varphi)=\left\{f\in L^2(\R): \{f(k)\}_{k\in\Z}\in \ell^2(\Z) ,\, f(x)=\sum_{k\in \Z} f(k)
s(x-k)\right\}.$$
\end{remark}

\bigskip

\begin{proof}[Proof of Theorem~\ref{t1}]

Let $f\in\mathcal{U}$ belong to a sampling space $V(\varphi)$ with
$\varphi\in\mathcal{U},$ and let $s$ be its sampling function.

$a)$ is true by Lemma~\ref{lemma}, and hence $Z_f(0,\omega)$ is
well defined, moreover, it is in $L^2[0,1].$

On the other side we have
$$G_f(\omega)=|Z_f(0,\omega)|^2G_s(\omega)\,\,\,\,\text{
a.e.} \, \omega \in \R,$$

hence, we can write (by Lemma~\ref{L3} $iv)$)

$$\frac{G_f(\omega)} {|Z_f(0,\omega)|^2}=
  G_{s}(\omega) \,\,\,\,\text{
a.e.} \, \omega \in E_f,$$ and since $E_f \subseteq E_s$ (up to a
set of measure zero), the right hand side of the  preceding
equation is bounded above and below
by $A$ and $B$ a.e. $\omega\in {E_f}$, and $b)$ follows.\\

For $c),$ note that

$$\int_{E_f} \frac{| \hat{f}(\omega)|}{|Z_f(0,\omega)|} d \omega=
\int_{E_f} |\hat{s}(\omega)| d \omega \leq\int_{-\infty}^{+\infty}
|\hat{s}(\omega)| d \omega,$$

so it suffices to show that $s\in\mathcal{U}.$ We have

$$\hat{\varphi}(\omega)=Z_{\varphi}(0,\omega)\hat{s}(\omega).$$

By Lemma~\ref{L3} $i),\,Z_{\varphi}(0,\cdot)$ is bounded below
(and above) a.e. in $E_{s}.$ Hence
$$\int_{\R}|\hat{s}(\omega)| d\omega=\int_{\R\cap E_{s}} |\hat{s}(\omega)| d\omega
\leq C \int_{\R\cap E_{s}}|\hat{\varphi}(\omega)|
d\omega=C\int_{\R}|\hat{\varphi}(\omega)| d\omega,$$

which proves our claim. We will now prove $d).$

By the Fourier Inversion Formula, for all $x \in \R$ and all $k\in
\Z$

$$s(x+k)=\int_{- \infty}^{+\infty} \hat{s}(\omega) e^{2 \pi i\omega (x+k)} d\omega=
\int_0^1 \sum_{j\in\Z} \hat{s}(\omega +j)e^{2\pi i(\omega+j)x}e^{2\pi i \omega k}
 d \omega.$$

Let $h(\omega,x)=\sum_{j \in \Z} \hat{s}(\omega + j)e^{2 \pi
i(\omega+j)x}.$ Using  Proposition \ref{SZ99} and Lemma \ref{L1}, $\forall x \in \R,$

$$\sum_{k\in\Z}\left|\int_0^1 h(\omega,x) e^{2\pi i \omega k} d
\omega \right|^2=\sum_{k \in \Z} |s(x+k)|^2 < L < +\infty. $$

On the other hand,
$$ \int_0^1 |h(\omega,x)| d \omega \leq \int_{-\infty}^{+\infty}
|\hat{s}(\omega)| d \omega,$$
so $h(\omega,x)\in L^1[0,1]$ as a function of $\omega.$\\
Then for all
$x\in \R,\, h(\cdot,x) \in L^2[0,1],$ and using
Parseval's equality, $\|h(\cdot,x)\|_2^2< L.$\\
So, since $\hat{s}(\omega)=\frac{\hat{f}(\omega)}{Z_f(0,\omega)}$
for $\omega\in E_f,$ we have
$$h(\omega,x)\chi_{[0,1]\cap E_f}(\omega)=e^{2\pi i \omega x} \frac{\sum_{s \in \Z}
\hat{f}(\omega + s)e^{2
\pi i sx}}{Z_f(0,\omega)}\chi_{[0,1]\cap E_f}(\omega),$$
and $d)$ follows.\\
This completes the proof of one implication. We will prove now the
sufficient conditions.
\bigskip
Define
\begin{equation}\label{3}
\hat{s}(\omega)=
\begin{cases}
\frac{\hat{f}(\omega )}{Z_f(0,\omega)}&
\text{for } \omega \in E_f \\
1 & \text{for }\omega\in [0,1]\setminus E_f\\
0 & \text{otherwise. }
\end{cases}
\end{equation}
Observe that if
$\omega \in E_f$, then by condition $b)\, |Z_f(0,\omega )|>0.$

We will prove that $V(s)$ is a sampling space with sampling
function $s$ and that $f\in V(s).$

First we will see that $\hat{s}\in L^2(\R)$ and that
$\{s(\cdot-k)\}_{k\in\Z}$ is a frame sequence.\\
For almost every $\omega \in E_f,$
$$\sum_{j\in \Z} |\hat{s}(\omega + j)|^2=
\sum_{j\in \Z}\frac{ |\hat{f}(\omega + j)|^2}{|Z_f(0,\omega)|^2}
=\frac{\sum_{j\in \Z} |\hat{f}(\omega +
j)|^2}{| Z_f(0,\omega)|^2},$$

therefore, using $b)$, it follows  that
\begin{equation}\label{5}
A \leq\sum_{k\in \Z} |\hat{s}(\omega + k)|^2 \leq B \,\,\,\text{
a.e. } \omega \in E_f.
\end{equation}
Now, since by definition of $s,$ $\sum_{k \in \Z}
|\hat{s}(x+k)|^2=1$ in $\R\setminus E_f,$ we have that

$$\min\{1,A\} \leq\sum_{k\in \Z} |\hat{s}(\omega + k)|^2 \leq
\max\{1,B\} \,\,\,\text{ a.e. } \omega \in \R,$$

and so $\{s(\cdot-k)\}_{k\in\Z}$ is a frame sequence. (Moreover,
this shows that it is a Riesz sequence).

Next we show that $\hat{s}\in L^1(\R)$:

$$\int_{- \infty}^{+ \infty} |\hat{s}(\omega)| d\omega=
\int_{E_f} \frac{| \hat{f}(\omega)|}{|Z_f(0,\omega)|} d\omega +
\int_{[0,1]\setminus E_f}d\omega.$$

So by $c)$,  $\hat{s} \in L^1(\R)$ and then $s$ is continuous.\\
Now we will see that $f \in V(s).$

$$\hat{f}(\omega)=Z_f(0,\omega) \frac{\hat{f}(\omega )}
{Z_f(0,\omega)}=Z_f(0,\omega) \hat{s}(\omega) \;\; \text{ a.e. }
\; \omega \in E_f,$$ hence using $b)$ we obtain

$$\hat{f}(\omega)=Z_f(0,\omega)\hat{s}(\omega) \;\;
\text{ a.e. }  \; \omega \in \R,$$ and therefore,
$$f(x)=\sum_{k \in \Z}f(k)s(x-k).$$

To complete the proof, it is enough to verify that the function
$s$ satisfies $ii)$ of Proposition~\ref{SZ99}. To see that
$\sum_{k \in \Z} |s(x+k)|^2$ is uniformly bounded, observe that
for all $x\in \R$ and all $k\in \Z$ we have
\begin{align*}
s(x+k) & =  \int_{- \infty}^{+\infty} \hat{s}(\omega) e^{2 \pi i\omega (x+k)
} d\omega \\*[3mm]
       &  =  \int_{ E_f} \frac{ \hat{f}(\omega )}{Z_f(0,\omega)}
 e^{2 \pi i\omega(x+k)} d \omega +\int_{[0,1]\setminus E_f}
 e^{2\pi i\omega(x+k)} d\omega \\*[3mm]
        & =  \int_{ [0,1] \cap E_f} \frac{\sum_{j \in \Z} \hat{f}(\omega +
j)e^{2 \pi i(\omega+j)(x+k)}}{Z_f(0,\omega)}  d \omega +
\int_{[0,1]\setminus E_f} e^{2\pi i\omega(x+k)} d\omega.
\end{align*}

Let us call $$g(\omega,x)=\frac{\sum_{j \in \Z} \hat{f}(\omega +
j)e^{2 \pi i(\omega+j)x}}{Z_f(0,\omega)}\chi_{[0,1]\cap
E_f}(\omega)+ e^{2\pi i\omega x}\chi_{[0,1]\setminus E_f}.$$ Then

$$s(x+k)=\int _0^1 g(\omega,x)e^{2\pi i k \omega} d \omega.$$
 Note
that, by $d)$, for all $x\in \R,\, g(\cdot,x) \in L^2[0,1]$ and
$\|g(\cdot,x)\|_2^2<L+1.$

By
Parseval's equality we have,

$$\sum_{k \in \Z} |s(x+k)|^2= \|g(\cdot,x)\|^2_2 < L+1 \,\,\,\forall x \in \R.$$

Since $\{f(k)\}_{k\in \Z}\in \ell^2(\Z) $ and $\hat{f}\in
L^1(\R)$, the Poisson Summation Formula holds for $f,$ that is:
$$Z_f(0,\omega)=\sum_{k\in\Z}\hat{f}(\omega +k)\,\,\,\text{ a.e. } \omega \in [0,1].$$

Therefore,

$$s(k)=\int_{-\infty}^{+\infty}\hat{s}(\omega)e^{2\pi ik \omega} d\omega=
\int_{E_f}\frac{\hat{f}(\omega)}{Z_f(0,\omega)}e^{2\pi ik \omega} d\omega
+\int_{[0,1]\setminus E_f}e^{2\pi i k\omega} d\omega$$

$$=\int_{[0,1]}e^{2\pi ik \omega} d\omega=\delta_{0k}, $$

hence $Z_{s}(0,\omega)=1.$\\

So using Proposition~\ref{SZ99}, it follows that $V(s)$ is a
sampling space. Moreover, by $iv)$ of Lemma~\ref{L3}, $s$ is its
sampling function.
\end{proof}
Note that for the proof of the necessary condition of
Theorem~\ref{t1} we did not use that $f\in\mathcal{U},$ hence
Proposition~\ref{nece} is immediate.

\bigskip
\section{Examples}

\begin{example}\label{Ej1}
It is easy to see that the $\sigma$ band-limited functions, for
which the Shannon Sampling Theorem holds, satisfy the conditions  of Theorem~\ref{t1}.\\
For instance, let $\sigma=\frac{1}{2}$. Then $\sum_{k\in\Z}
\hat{f}(\omega+k)=\hat{f}(\omega)\in L^2[0,1],$ so
$\{f(k)\}_{k\in\Z}\in \ell^2(\Z).$ Then, using Poisson Summation
Formula, $Z_f(0,\omega)=\hat{f}(\omega),$
and conditions $b), c) $ and $d)$ follow.\\

\end{example}

In \cite{SZ04} the following sufficient conditions for a function
to belong to sampling space are given:

\begin{proposition}\label{sufc}(\cite{SZ04})
Let $f\in\mathcal{U}$. If there are positive constants $A$ and $B$
such that

$$ A\left|\sum_{k\in\Z}\hat{f}(\omega+k)\right|^2\leq
\sum_{k\in\Z}|\hat{f}(\omega+k)|^2, \,\,\, \text{ a.e.}$$ and
$$\left(\sum_{k\in\Z}|\hat{f}(\omega+k)|\right)^2\leq B \left|\sum_{k\in\Z}\hat{f}(\omega+k)\right|^2,\,\,\, \text{ a.e.}$$
then $f$ belongs to a sampling space.
\end{proposition}
Now we will show that the conditions of our Theorem~\ref{t1} are
weaker than these conditions.

\begin{example}\label{Ej2}
Define $f\in L^2(\R)$ such that:

$$\hat{f}(\omega)=\sum_{n=0}^{+\infty}\frac{(-1)^n}{n+1}
\chi_{[n,n+\frac{1}{2^n}]}(\omega).$$

Then
$$\sum_{k\in \Z}\hat{f}(k)=\sum_{n=0}^{+\infty}\frac{(-1)^n}{n+1}$$
and for every $\omega\in (1/2^{n+1},1/2^n]:$
$$\sum_{k\in \Z} \hat{f}(\omega +
k)=\sum_{k=0}^{n}\frac{(-1)^k}{k+1},\,\,\,\,\,\sum_{k\in
\Z} |\hat{f}(\omega + k)|=\sum_{k=0}^{n}\left|\frac{1}{k+1}\right|$$
$$\sum_{k\in \Z} |\hat{f}(\omega +
k)|^2=\sum_{k=0}^{n}\left|\frac{1}{k+1}\right|^2.$$

 It is easy to
see that $\hat{f}\in L^1(\R)$ and that $\sum_{k\in
\Z}\hat{f}(\omega+k)\in L^2[0,1].$ So $\{f(k)\}_{k\in\Z}\in
\ell^2(\Z),$ and using Poisson Summation Formula it follows that
$f$
 satisfies the conditions $b), c)$ and $d)$ of Theorem~\ref{t1}.\\
 On the other side, $f$ doesn't satisfy
 the conditions of Proposition~\ref{sufc}.
\end{example}

\begin{remark}
We need to assume that $f,\varphi\in \mathcal{U}$ since in our proof we use the
Fourier Inversion Formula and the Poisson Summation Formula.
\end{remark}

Note that there exist functions belonging to sampling spaces  that  are not in
$\mathcal{U}$, as we will see in the following example:
\begin{example}
Let
$$s(x)= \begin{cases} -\sin(\pi x) &\text{for }\,\,x\in [-1,-1/2]\\
                         1& \text{for }\,\,x \in [-1/2,1/2]\\
                         \sin(\pi x )&\text{for }\,\,x\in [1/2,1]\\
                         0&\text{for }\,\,x \notin  [-1,1]
\end{cases}$$
This is an example of a sampling function that doesn't belong to
$\mathcal{U}.$
  Applying Lemma 7.3.3 and Corollary 7.3.4 of \cite{Chr03},
 $\{s(\cdot - k), k\in \Z\}$
 is a Riesz basis for the space $V(s).$
 Since $s$ is a compactly supported function,\\
 $$\sum_{k\in \Z} |s(x-k)|^2 <L<+\infty \,\,\,\forall x \in \R,
\text{ and }\sum_{k \in \Z} s(k) e^{-2 \pi i k\omega}=s(0)=1.$$
Using Proposition~\ref{SZ99}, $V(s)$ is a sampling space with sampling
function $s.$ But $\hat{s} \notin L^1(\R),$ moreover, for every
$f\in V(s)$, $\hat{f}\notin L^1(\R).$
\end{example}

\vspace{13pt}

\centerline{ACKNOWLEDGMENT}

\vspace{13pt}

The research of the authors is partially supported by Grants
PICT03-15033 and UBACYT X058.

The authors want to thank Ursula Molter for reading carefully the
manuscript and for her suggestions that helped to improve the
paper. Part of this research was carried out while C. Cabrelli and
S. Heineken where visiting NuHAG, University of Vienna, and the
Erwin Schr\"odinger International Institute for Mathematical
Physics, Vienna, during the special semester on Modern Methods of
Time-Frequency Analysis 2005. They want to thank the organizers
for their generous hospitality.

\end{document}